\documentclass[12pt,a4paper]{article}
\usepackage{amsmath}
\usepackage{amssymb}

\def\beq{\begin{equation}}
\def\eeq{\end{equation}}

\def\be{\begin{displaymath}}
\def\ee{\end{displaymath}}

\newtheorem{thm}{Theorem}
\newtheorem{col}{Corollary}
\newtheorem{lem}{Lemma}

\title{
On a Functional Differential Equation of Determinantal Type
}

\author{H. W. Braden\thanks{E-mail:hwb@ed.ac.uk}\ \ and 
J.G.B. Byatt-Smith\thanks{E-mail:Byatt@ed.ac.uk}\\
\normalsize
\em Department of Mathematics and Statistics,\\
\normalsize
\em The University of Edinburgh, \\
\normalsize
\em Edinburgh, UK \\
}

\begin{document}

\renewcommand{\thepage}{}
\begin{titlepage}

\maketitle
\vskip-9.5cm
\hskip10.4cm
MS-98-005
\vskip8.8cm

\begin{abstract}
We solve the functional equations
$$
\begin{vmatrix}
1              & 1              & 1              \\
f(x)           & f(y)           & f(z)           \\
f\sp{\prime}(x)& f\sp{\prime}(y)& f\sp{\prime}(z)\\
\end{vmatrix}
=0,\quad\quad
\begin{vmatrix}
1              & 1              & 1              \\
f(x)           & g(y)           & h(z)           \\
f\sp{\prime}(x)& g\sp{\prime}(y)& h\sp{\prime}(z)\\
\end{vmatrix}
=0,
$$
for suitable functions $f$, $g$ and $h$ subject to $x+y+z=0$.
These equations essentially characterise the Weierstrass $\wp$-function
and its degenerations.

\end{abstract}

\begin{flushleft}
\end{flushleft}

\vfill
\end{titlepage}
\renewcommand{\thepage}{\arabic{page}}

\section{Introduction}
The purpose of the following note is to present a simple and direct proof of
\begin{thm}
Let f be a three-times differentiable function satisfying the functional
equation
\begin{equation}
\begin{vmatrix}
1              & 1              & 1              \\
f(x)           & f(y)           & f(z)           \\
f\sp{\prime}(x)& f\sp{\prime}(y)& f\sp{\prime}(z)\\
\end{vmatrix}
=0,\quad\quad  x+y+z=0.
\label{detdiff}
\end{equation}
Up to the manifest invariance
\begin{equation*}
f(x)\rightarrow \alpha f(\delta x)+\beta,
\label{equiv}
\end{equation*}
the solutions of (\ref{detdiff}) are one of $f(x)=\wp(x+d)$,
$f(x)=e\sp{x}$ or $f(x)=x$.
Here $\wp$ is the Weierstrass $\wp$-function and 
 $3 d$ is a lattice point of the $\wp$-function.
\end{thm}
In fact our approach gives a simple proof of
\begin{thm}
Let $f$, $g$ and $h$ be  three-times differentiable functions 
satisfying the functional equation
\begin{equation}
\begin{vmatrix}
1              & 1              & 1              \\
f(x)           & g(y)           & h(z)           \\
f\sp{\prime}(x)& g\sp{\prime}(y)& h\sp{\prime}(z)\\
\end{vmatrix}
=0,\quad\quad x+y+z=0.
\label{reddetdiff}
\end{equation}
Up to the manifest invariance
\begin{equation*}
f(x)\rightarrow \alpha f(\delta x+\gamma_1)+\beta,\ \,
g(x)\rightarrow \alpha g(\delta x+\gamma_2)+\beta,\ \,
h(x)\rightarrow \alpha h(\delta x+\gamma_3)+\beta,
\end{equation*}
where $\gamma_1+\gamma_2+\gamma_3=0$, the nonconstant solutions of
(\ref{reddetdiff}) are given by $f(x)=g(x)=h(x)=e\sp{x}$, $x$, or $\wp(x)$.
If (say) $h(z)$ is a constant then either
\begin{enumerate}
\item One of the functions $f(x)$ or $g(y)$ is the same constant as
$h(z)$, in which case the remaining function is arbitrary, or
\item $f(x)=g(x)=e\sp{x}$.
\end{enumerate}

\end{thm}
\noindent{\bf Remarks:} $(i)$
In fact the exponential and linear function solutions satisfy 
(\ref{detdiff}) and (\ref{reddetdiff}) without the constraint
$x+y+z=0$.\newline
$(ii)$ The theorems immediately give the general analytic
solutions to the same functional equations viewed as functions of a 
complex variable, showing that the solutions are in fact
meromorphic.\newline
$(iii)$ The  arbitrary constant $\delta$ in the invariance of (\ref{equiv})
is accommodated in the 
Weierstrass $\wp$-function solution  by the homogeneity relation
$\wp(tx;t\sp{-4}g_2,t\sp{-6}g_3)=t\sp{-2}\wp(x;g_2,g_3)$.
\newline
$(iv)$ Weierstrass has shown \cite{PS} that any meromorphic function
possessing an algebraic addition formula is either an elliptic function
or is of the form $R(z)$ or $R(e\sp{\lambda z})$, where $R$ is a rational
function. A priori the functional equation (\ref{detdiff}) is distinct
from assuming $f$ possesses an algebraic addition formula.

Interest in (\ref{detdiff}) arises from a question of mathematical
physics: 
What one-dimensional quantum mechanical models with pair-wise
interactions have ground states of product type? 
In addressing this question \cite{suth, cal} the functional equation 
\begin{equation}
F(x)F(y) +F(y) F(z)+F(z) F(x) =G(x)+G(y)+G(z),
\label{factfun}
\end{equation}
appears (with $x+y+z=0$), where on physical grounds $F(x)$ is taken to be odd.
By applying $(\partial_x-\partial_y)\, \partial_x\,\partial_y$ to this
equation we obtain (\ref{detdiff}) with $f(x)=F\sp\prime(x)$
and similarly (\ref{detdiff}) may be integrated to yield (\ref{factfun}).
The models that arise in the solution of this question include
the Calogero-Moser-Sutherland models. They are rich in 
interesting mathematics involving representation theory, harmonic
analysis and special functions (see for example \cite{heckman}); the
classical analogues of the models yield completely integrable Hamiltonian
systems. Now the solutions that yield these models were obtained assuming the
function $F(x)$, or equivalently $f(x)$, to be meromorphic with at 
least one pole. With such an assumption it is easy to show that $f(x)=
\wp(x +d)$ is the general solution, and further requiring $f(x)$ to be
even dictates $d=0$. Although one can in fact show that there are no 
(nonconstant)
even entire solutions to (\ref{detdiff}), so answering the physical question,
the general solution to (\ref{detdiff}) appeared difficult to obtain.
Indeed, the functional equation (\ref{reddetdiff}) was introduced in
\cite{bp} as a means to understand (\ref{detdiff}),\footnote{We are grateful to
V.M. Buchstaber for informing us on this matter.} but the proof of
the main theorem of \cite{bp} fails to result in a direct proof of Theorem 1.
Here we present a simple and direct proof of (\ref{detdiff}) that
also allows a new and simpler proof of (\ref{reddetdiff}).

\section{Proofs}
The strategy of our proof is first to isolate a necessary condition for
nonconstant solutions of (\ref{detdiff}) and (\ref{reddetdiff}) in the
form of a differential equation. The solutions of this
differential equation  are given in terms of the Weierstrass $\wp$-function 
or one of its degenerations, and the second step is to determine the various 
free parameters that arise in the solution. This will prove the theorems
for nonconstant functions $f$, $g$ and $h$. Finally we consider the case
where some of these functions are constant. 

While the approach to solving a functional equation via an associated 
differential equation(s) is standard \cite{Acza}, the simplicity
of our proof depends on one rather nonobvious step that we wish to
highlight in advance.
First, by taking (\ref{reddetdiff}) and various of its derivatives  we 
obtain an equation of the form $F(f(x),g(y))=0$ involving $f(x)$, $g(y)$
and their derivatives. Viewing $x$ (say $x=x_0$) as fixed we have in general a
nonlinear ordinary differential equation
for $g(y)$ which may in principle be solved.
While such a solution satisfies $F(f(x_0),g(y))=0$, it needn't in general
satisfy $F(f(x),g(y))=0$, the various derivatives 
$\partial_x\sp{n}\, F(f(x),g(y))|_{x_0}=0$ 
(supposing they exist) yielding further differential equations for $g(y)$.
These further equations give us restrictions on the  allowed functions
$g(y)$ and between them one may eliminate various derivatives of $g$
appearing.
For example, given $F(f(x_0),g(y))=0$ and $\partial_x\, F(f(x),g(y))|_{x_0}=0$
one could could choose to eliminate the highest derivative of $g$ appearing;
similarly one can use  further partial derivatives to eliminate
additional derivatives of $g$. 
If we suppose that $F(f(x),g(y))=0$ determines $g(y)$ then
the equation $F(f(x_0),g(y))=0$, together with an appropriate number of further
derivatives,  also determines $g(y)$. 
The nonobvious step in our proof is to provide a functional
$F(f(x),g(y))=0$ that alone readily gives $g(y)$. We will
remark in the course of the proof when this is done, and note the
several perhaps surprising simplifications that follow.

\begin{lem} Let $f$, $g$ and $h$ be three-times differentiable, 
nonconstant functions that satisfy (\ref{reddetdiff}).
Then  each satisfies the (same) differential equation
\begin{equation}
w'(x)^2=p_3\, w(x)^3+p_2\, w(x)^2+p_1\, w(x)+p_0.
\label{diffsola}
\end{equation}

\end{lem}
\noindent{\it Proof:} 
We begin by deriving several algebraic consequences of
assuming that the (nonconstant) functions $f$, $g$ and $h$ of 
(\ref{reddetdiff}) are $N$-times differentiable.
Let $1\le k,l,s\le N$. 
The algebraic identities we obtain yield a large
supply of functional equations involving only the functions $f(x)$, $g(y)$
and their derivatives.
We will obtain (\ref{diffsola}) by eliminating an appropriate derivative.
A minimum choice of $N=3$ will arise in the proof. 

Set $\partial=\partial_y -\partial_x$ and let
\begin{equation}
a_k=\partial\sp{k-1} ( g(y)-f(x) ), \quad
b_k=\partial\sp{k}   ( g(y)+f(x) ), \quad
c_k=\partial\sp{k}   ( g(y)\, f(x)).
\end{equation} 
Then differentiation of (\ref{reddetdiff}) yields $N$ equations,
\begin{equation*}
a_k\, h\sp{\prime}(z)-b_k\, h(z) +c_k=0.
\end{equation*}
Comparing any two of these equations shows
\begin{equation}
(a_k \, b_l - a_l \, b_k ) h(z) = a_k \, c_l - a_l \, c_k  ,
\quad\quad
(a_k \, b_l - a_l \, b_k )h\sp{\prime}(z)= b_k \, c_l - b_l \, c_k ,
\label{twoeq}
\end{equation}
while comparison of any three yields
\begin{equation}
\begin{vmatrix}
a_k   & a_l  & a_s              \\
b_k   & b_l  & b_s         \\
c_k   & c_l  & c_s         \\
\end{vmatrix}
=0.
\label{threedet}
\end{equation}
Consider $z=-x-y$ as a function of $x$ and $y$.
Differentiating the first of equations (\ref{twoeq}) with respect to
$y$  say, and comparing with the second  equation results in 
\begin{equation}
\begin{vmatrix}
a_k   & a_l  & a_k \, \partial_y a_l -a_l \, \partial_y a_k \\
b_k   & b_l  & a_k \, \partial_y b_l -a_l \, \partial_y b_k         \\
c_k   & c_l  & a_k \, \partial_y c_l -a_l \, \partial_y c_k
+b_k \, c_l - b_l \, c_k  \\
\end{vmatrix}
=0.
\label{addet}
\end{equation}
Similar expressions result upon differentiation with respect $x$.

Observe that at this stage (\ref{threedet}) and (\ref{addet})  and
their linear combinations provide us with many
functional equations involving only the functions $f(x)$, $g(y)$
and their derivatives. 
In particular we may eliminate various combinations of derivatives between
them. For example, by taking $k=1$, $l=2$ and $s=3$ the linear
combination
(\ref{addet}) -- $a_1\times$(\ref{threedet}),
\begin{equation}
\begin{vmatrix}
a_1   & a_2  & a_1\, g''(y)-a_2\, g'(y)-a_1\, a_3 \\
b_1   & b_2  & a_1\, g'''(y)-a_2\, g''(y)-a_1\, b_3 \\
c_1   & c_2  & a_1 \, \partial_y c_2 -a_2 \, \partial_y c_1 
-a_1\, c_3 +(b_1\, c_2- b_2\, c_1 )        \\
\end{vmatrix}
=0,
\end{equation}
cancels each of the  $g'''(y)$ derivatives appearing in the third column.
This expression, quadratic in $g''(y)$, factorises to give
\begin{equation}
\begin{split}
0=& \left( f'(x)\,g'(y) - {{g'(y)}^2} - f(x)\,g''(y) + g(y)\,g''(y)
        \right)\times \\
&\bigg( 3{{f'(x)}^3} - 3f'(x){{g'(y)}^2} - 
       4f(x)f'(x)f''(x) + 4g(y)f'(x)f''(x)  + {{g(y)}^2}f^{(3)}(x)\\
&- 2f(x)f'(x)g''(y)  + 2g(y)f'(x)g''(y) +  
{{f(x)}^2}\,f^{(3)}(x)  - 2\,f(x)\,g(y)\,f^{(3)}(x)\bigg).
\end{split}
\label{eqfact}
\end{equation}
This elimination of the  $g'''(y)$ derivatives gives us the
equation $F(f(x),g(y))=0$ we choose to work with.
The factorisation we encounter is one of the simplifications we drew 
attention to earlier.  
Of course we could have taken the  $k=1$, $l=2$ and $s=3$ form of
(\ref{addet}) as our equation $F(f(x),g(y))=0$. 
It appears that this, together with one further partial derivative,
$\partial_x\, F(f(x),g(y))=0$ is sufficient to determine $g(y)$, but at the 
expense of far greater work. In particular there is no similar
factorisation  to that we encountered above. 
The proof presented was devised to circumvent the difficulties of this latter
route.

Now the first term appearing on the right of (\ref{eqfact}) may be written as
\begin{equation*}
{{\left( f(x) - g(y) \right) }^2}\,
\frac{d}{dy}\bigg( {\frac{f'(x) - g'(y)}{f(x) - g(y)}}\bigg), 
\end{equation*}
while the second term may be expressed as 
\begin{equation*}
{\frac{{{\left( f(x) - g(y) \right) }^4}}{g'(y)}}
\frac{d}{dy}\bigg(
{\frac{f'(x)\,\left( {{f'(x)}^2} - {{g'(y)}^2} \right) }
     {{{\left( f(x) - g(y) \right) }^3}}} - 
   {\frac{2\,f'(x)\,f''(x)}{{{\left( f(x) - g(y) \right) }^2}}} +
   {\frac{f^{(3)}(x)}{f(x) - g(y)}}
\bigg).
\end{equation*}
Our assumption of nonconstant solutions means that $f(x)\not\equiv g(y)$ and
$g'(y)\not\equiv0$. Thus the vanishing of (\ref{eqfact}) means either
\begin{equation*}
{\frac{f'(x) - g'(y)}{f(x) - g(y)}}=C_1(x)
\end{equation*}
or
\begin{equation*}
{\frac{f'(x)\,\left( {{f'(x)}^2} - {{g'(y)}^2} \right) }
     {{{\left( f(x) - g(y) \right) }^3}}} -
   {\frac{2\,f'(x)\,f''(x)}{{{\left( f(x) - g(y) \right) }^2}}} +
   {\frac{f^{(3)}(x)}{f(x) - g(y)}}
=C_2(x),
\end{equation*}
according to whether the first or second terms in (\ref{eqfact}) vanish.
Therefore, assuming only  third derivatives exist, we have shown (after 
rearranging) that either
\begin{equation}
g'(y)=l_1\, g(y)+l_0, \quad{\rm or}\quad
g'(y)^2=p_3\, g(y)^3+p_2\, g(y)^2+p_1\, g(y)+p_0.
\label{diffsol}
\end{equation}
Both are  cases of (\ref{diffsola}).
A consequence of (\ref{diffsol}) is that derivatives to all orders exist
for solutions of either of these differential equations. 
The (analytic) solutions of these differential equations may be expressed
(generically) in terms of the exponential and Weierstrass $\wp$-function,
 and indeed
the exponential solution of the linear differential equation corresponds to a 
degeneration of the $\wp$-function differential equation.
Further, the identical argument but upon interchanging the role
of $y$ and $x$ shows that $f(x)$ is also a $\wp$-function or a degeneration.
We have however yet to establish that $f$ and $g$ satisfy the same 
differential equation. Before turning to this there is one point that
needs to be clarified.

In principle it is possible for
the function $g(y)$ giving the vanishing of (\ref{eqfact}) to be a solution of
one of the differential equations (\ref{diffsol}) in one domain and satisfy
the other differential equation outside of it, with $g(y)$ and its first
three derivatives matching at any boundary. Such a possibility does not
arise in our problem. One can readily show that  matching  a solution
$g_1(y)$ of the first differential equation (\ref{diffsol})  with a
solution $g_2(y)$ of the second differential equation at a point
$y_0$, and requiring the first four derivatives to agree at this point,
entails $g_1(y)\equiv g_2(y)$. Thus  the solutions to the vanishing
of (\ref{eqfact}) being envisaged perforce have discontinuous fourth
derivative at such points $y_0$.
Now for our problem, the coefficients $l_i$ and $p_i$
of the differential equations (\ref{diffsol}) are in fact functions of
$x$, given explicitly below. We will shortly see that it is this
aspect of our problem that rules out
the functions $g(y)$ envisaged in this paragraph.

We now establish that the functions $f$ and $g$ satisfy the same 
differential equation by determining the coefficients $l_i$ and $p_i$
appearing above in two different ways. By directly differentiating
(\ref{diffsol}) one obtains (for example) that
\begin{equation*}
l_1= {\frac{g''(y)}{g'(y)}},\qquad
p_3= {\frac{ g'(y)\,g^{(4)}(y)-g''(y)\,g^{(3)}(y)}{3\,{{g'(y)}^2}}}.
\end{equation*}
(The remaining coefficients will be listed below.)
Alternately the coefficients may be determined by evaluating the
functions $C_{1,2}(x)$ arising upon integration.
These functions  may be determined in a variety of ways. If (\ref{eqfact})
vanishes, so does the derivative with respect to $x$ of the right-hand side
of the equation. Substituting either of the expressions for $g''(y)$
corresponding to the two terms of (\ref{eqfact}) into this derivative
equation yield equations for $g'(y)$. Comparison with the equations
defining $C_{1,2}(x)$ show (with a little work)
\begin{equation*}
C_1(x)= {\frac{f''(x)}{f'(x)}},\qquad
C_2(x)= {\frac{ f'(x)\,f^{(4)}(x)-f''(x)\,f^{(3)}(x)}{3\,{{f'(x)}^2}}}.
\end{equation*}
(The same expressions arise by applying
L'Hospital's rule to the equations defining $C_{1,2}(x)$.)  
The coefficients $l_i$ and $p_i$ are then determined to be
\begin{equation}
\begin{split}
l_0&=\frac{f'(x)^2 -f(x)\,f''(x)}{f'(x)}\\
l_1&={\frac{f''(x)}{f'(x)}}\\
p_0&=f'(x)^2 - 2 f(x)\,f''(x) +\frac{f(x)^2 f^{(3)}(x)}{f'(x)}
-\frac{f(x)^3}{3\, f'(x)^3}\left( f'(x)\,f^{(4)}(x)- f''(x)\,f^{(3)}(x)\right)
\\
p_1&= 
2\,f''(x) - {\frac{2\,f(x)\,f^{(3)}(x)}{f'(x)}} 
+\frac{f(x)^2}{f'(x)^3}\left(f'(x)\,f^{(4)}(x)- f''(x)\,f^{(3)}(x)\right) 
\\
p_2&=\frac{f^{(3)}(x)}{ f'(x)} -\frac{f(x)}{f'(x)^3}
\left(f'(x)\,f^{(4)}(x)- f''(x)\,f^{(3)}(x)\right)
\\
p_3&={\frac{ f'(x)\,f^{(4)}(x)-f''(x)\,f^{(3)}(x)}{3\,{{f'(x)}^3}}}
\\
\end{split}
\end{equation}
Thus, for example, we have shown that
\begin{equation*}
l_1= {\frac{g''(y)}{g'(y)}}={\frac{f''(x)}{f'(x)}},
\end{equation*}
and so $l_1$ is in fact a constant. Similarly $p_3$ is seen above to be
a constant. Further, these constants are determined by expressions
symmetric in the interchange of $f$ and $g$ and the same is true of the
remaining coefficients. This symmetry means that $f$ and $g$ will
satisfy the same differential equation. Also, if we had chosen to
work with $f$ and $h$ the above argument shows that they satisfy the
same differential equation. Therefore, each of $f$, $g$ and $h$ satisfy
the same differential equation. Finally observe that the coefficient
$p_3$, which is constant, involves the fourth derivative of the
functions satisfying (\ref{diffsol}) for arbitrary position. This
expression shows one cannot
construct a constant $p_3$ from two (nonconstant) solutions of (\ref{diffsol}) 
whose first three derivatives agree and whose
fourth derivatives differ. This rules out having solutions to
(\ref{diffsol}) that satisfy the first equation on one region and the
second equation elsewhere, and having discontinuous fourth derivative on the
boundaries.

\hfill$\square$
\begin{col}
A nonconstant solution $f(x)$  of (\ref{detdiff}) satisfies the differential 
equation
(\ref{diffsola}).
\end{col}
\noindent{\bf Remarks:}$(i)$ If one was happy to assume the five-times
differentiability of the function $f(x)$ satisfying (\ref{detdiff}) then
this result may be obtained very quickly. Let $f=g$ in the lemma. We need
proceed no further than (\ref{addet}). Simplification arises because we
have only the one function  and its derivatives appearing throughout and
we may further let $y=x$ yielding (nonlinear) differential equations
satisfied by $f(x)$. Thus, for example,
taking $k=2$, $l=4$ and upon setting $y=x$ in (\ref{addet}) we obtain 
\begin{equation}
0=
(f')\sp6 \,
\left( \frac{f''}{f'}\right )' \,
\left ( \frac{1}{f'}\left( \frac{f'''}{f'} \right )' \right )' .
\label{eqf}
\end{equation}
(This in fact corresponds to first nonzero term appearing in a Taylor
series expansion around $y=x+\epsilon$ in the $k=1$, $l=2$ expressions of the
lemma.)
The two final two factors here correspond to the two terms appearing in
(\ref{eqfact}); they may be straightforwardly integrated to yield the
corollary. \newline
$(ii)$
Another (again not straightforward) route to solving (\ref{detdiff})
is first to express this functional equation in a rather different manner.
Parameterise $x+y+z=0$ by $x=\xi +\eta$, $y=\xi -\eta$ and $z=-2\xi$. Then
with $\Delta$ denoting the standard central difference operator (defined by
$\Delta  f(\xi)=f(\xi +\eta)-f(\xi -\eta)$) and $\mu$  denoting the
standard average operator (with  $\mu  f(\xi)=(f(\xi +\eta)+f(\xi -\eta))/2$)
then (\ref{detdiff}) may be written
$$
\Delta f(\xi)\, \mu f^{\prime}(\xi)-\Delta f^{\prime }(\xi)\, \mu f(\xi)=
\Delta f(\xi)\, F^{\prime }(z)-\Delta f^{\prime }(\xi)\, F(z).
$$
Here we have set $F(z)=f(-2\xi)$ and $F^{\prime
}(z)=-{\frac {d}{d\xi}F(z)}/2$ and we may view this as 
an equation in $\eta$ with $\xi$ fixed.
Expanding this equation as a series in $\eta$ leads to an infinite set of
relations which are first order in $F(z)$,
\begin{equation*}
a_k\, F\sp{\prime}(z)-b_k\, F(z) +c_k=0.
\end{equation*}
These equations are analogous to those involving $a_i$, $b_i$ and $c_i$ in 
the lemma. Here the  coefficients are
functions of $f(\xi)$ and its derivatives. The lowest order
coefficients are
\begin{equation*}
\begin{split}
a_1 &=f^{\prime},\quad b_1= f^{\prime \prime },\quad
c_1=ff^{\prime \prime }-f^{\prime 2}\\
a_2&=f^{\prime \prime \prime},\quad
b_2=f^{\prime \prime \prime \prime }, \quad
c_2=-4f^{\prime} f^{\prime \prime \prime}+ ff^{\prime \prime \prime
\prime }+3f^{\prime \prime 2}.
\end{split}
\end{equation*}
Elimination of $F(z)$ amongst these is also sufficient to obtain (\ref{eqf}).

At this stage we have found a necessary condition for nonconstant
functions
satisfying (\ref{reddetdiff}) and (\ref{detdiff}): the functions
must satisfy (\ref{diffsola}) and so \cite{WW} are of the form 
$$ \alpha \wp(\delta x+\gamma)+\beta,\qquad
 \alpha e\sp{\delta x}+\beta,\qquad
 \alpha x+\beta.$$
As remarked upon in the introduction, the exponential and linear solutions
clearly satisfy the functional equations without the restriction $x+y+z=0$
and we need only discuss the first solution here involving the $\wp$-function.
The noted invariances of the functional equations mean we need only
determine whether any restrictions must be placed on the translation parameter 
$\gamma$ unspecified by the differential equation.
\begin{lem} The function $f(x)= \wp(x+\gamma)$
satisfies (\ref{detdiff}) provided $3\gamma$ is a lattice point of the
$\wp$-function.

Similarly, the functions $f(x)=\wp( x+\gamma_1)$, $g(y)=\wp( y+\gamma_2)$,
$h(z)=\wp( z+\gamma_3)$ satisfy (\ref{reddetdiff}) provided 
$\gamma_1+\gamma_2+\gamma_3$ is a lattice point of the
$\wp$-function and these may be chosen so that 
$\gamma_1+\gamma_2+\gamma_3=0$.
\end{lem}
\noindent{\it Proof:}
The result follows from the following identity \cite{WW}
\begin{equation*}
\begin{vmatrix}
1              & 1              & 1              \\
\wp(a)           & \wp(b)           & \wp(c)           \\
\wp\sp{\prime}(a)& \wp\sp{\prime}(b)& \wp\sp{\prime}(c)\\
\end{vmatrix}
= 2 \frac{\sigma(a+b+c)\,\sigma(a-b)\,\sigma(b-c)\,\sigma(c-a)}{\sigma(a)^3\,
     \sigma(b)^3\,\sigma(c)^3},
\end{equation*}
where $\sigma$ is the Weierstrass sigma function that vanishes at the
lattice points of the $\wp$-function. Letting $a=x+\gamma$, $b=y+\gamma$, 
$c=z+\gamma$ we see that $f(x)= \wp(x+\gamma)$ 
satisfies (\ref{detdiff}) provided $3\gamma$ is a lattice point as stated.
Similarly we see the functions $f$, $g$ and $h$ given in the lemma
satisfy (\ref{reddetdiff}) provided $\gamma_1+\gamma_2+\gamma_3$ 
is again a lattice point. Suppose now $\gamma_1+\gamma_2+\gamma_3\equiv L$
is such a lattice point, and set $\gamma_3\sp\prime=-(\gamma_1+\gamma_2)$. 
The periodicity of the $\wp$-function means
$$ h(z)=\wp( z+\gamma_3)=\wp( z+\gamma_3-L)=\wp( z+\gamma_3\sp\prime)$$
and so we may choose the solutions of (\ref{reddetdiff}) so that
$\gamma_1+\gamma_2+\gamma_3=0$.
\hfill$\square$

It remains to discuss the case when at least one of the functions in
(\ref{reddetdiff}) is a constant which we may take to be $h(z)$.
\begin{lem}
Let $f$, $g$, $h$ satisfy (\ref{reddetdiff}). If $h(z)$ is a constant
then  either
\begin{enumerate}
\item One of the functions $f(x)$ or $g(y)$ is the same constant as
$h(z)$, in which case the remaining function is arbitrary, or
\item up to the invariance of (\ref{reddetdiff}), $f(x)=g(x)=e\sp{x}$.
\end{enumerate}
\end{lem}
\noindent{\it Proof:}
Using the invariance of (\ref{reddetdiff}) we may suppose without
loss of generality  that $h(z)=0$. Then (for all $x,y$)
\begin{equation}
0=
\begin{vmatrix}
f(x)           & g(y)          \\
f\sp{\prime}(x)& g\sp{\prime}(y)\\
\end{vmatrix}
=\left(\partial_y -\partial_x\right) f(x)\, g(y).
\label{constdiff}
\end{equation}
If neither of $f(x)$ or $g(y)$ vanish identically then 
$f(x)=e\sp{\delta x+\epsilon_1}$ and $g(y)=e\sp{\delta y+\epsilon_2}$,
with $\epsilon_{1,2}$ arbitrary and $\delta$ possibly zero (giving
the constant solutions). Using the invariance of (\ref{reddetdiff})
we may choose $f(x)=g(x)=e\sp{x}$.  Finally if (say) $f(x)=0$ 
we see that $g(y)$ is arbitrary.
\hfill$\square$

\section{Acknowledgements}
We wish to thank T. Potter, H. Ochiai and T. Oshima for correspondence
pertaining to these functional equations, V.M. Buchstaber for
various discussions and both
A. M. Davie and A. Olde Daalhuis for critically reading the manuscript.

\end{document}